\documentclass[11pt]{article}
\usepackage{jheppub}
\usepackage{bm,bbm}
\usepackage{amsmath}
\usepackage{mathtools}
\usepackage{graphics}
\usepackage{braket}
\usepackage{subcaption}

\usepackage{amsthm}

\edef\restoreparindent{\parindent=\the\parindent\relax}
\usepackage{parskip}
\restoreparindent
\usepackage{ascmac}
\usepackage{mathrsfs}
\usepackage{tikz}
\usetikzlibrary{patterns,intersections, calc,positioning,decorations.pathreplacing,decorations.pathmorphing,decorations.markings}
\tikzset{>=latex}
\tikzset{middlearrow/.style={
        decoration={markings,
            mark= at position 0.57 with {\arrow{#1}} ,
        },
        postaction={decorate}
    }
}
\def\d{{\rm d}}
\def\e{{\rm e}}

\def\i{{\rm i}}

\def\b0{\bm{0}_\perp}

\makeatletter
\newcommand*{\rom}[1]{\expandafter\@slowromancap\romannumeral #1@}
\makeatother

\DeclareFontFamily{U}{mathx}{\hyphenchar\font45}
\DeclareFontShape{U}{mathx}{m}{n}{
      <5> <6> <7> <8> <9> <10>
      <10.95> <12> <14.4> <17.28> <20.74> <24.88>
      mathx10
      }{}
\DeclareSymbolFont{mathx}{U}{mathx}{m}{n}
\DeclareFontSubstitution{U}{mathx}{m}{n}
\DeclareMathAccent{\widecheck}{0}{mathx}{"71}
\usepackage{accents}

\newtheorem{theorem}{Theorem}
\newtheorem{proposition}{Proposition}
\newtheorem{corollary}{Corollary}[theorem]

\title{On some identities for confluent hypergeometric functions and Bessel functions}

\author[a,b]{Yoshitaka Okuyama}

\affiliation[a]{
Department of Physics, Osaka University,\\
Machikaneyama-Cho 1-1, Toyonaka, Osaka 560-0043, Japan
}

\affiliation[b]{Department of Physics, Faculty of Science,
The University of Tokyo,\\
Bunkyo-Ku, Tokyo 113-0033, Japan}

\preprint{OU-HET-1167}

\abstract{
Mathematical functions, which often appear in mathematical analysis, are referred to as special functions and have been studied over hundreds of years. Many books and dictionaries are available that describe their properties and serve as a foundation of current science.
In this paper, we find a new integral representation of the Whittaker function of the first kind and show a relevant summation formula for Kummer's confluent hypergeometric functions. We also perform the specifications of our identities to link to known and new results.
}

\begin{document}
\maketitle

\section{Introduction}
The gamma function \cite[Chapter 5]{NIST23} is defined by the following integral representation for $\Re\,z>0$ \cite[(5.2.1)]{NIST23}:
\begin{align}\label{eq:gamma func def}
    \Gamma(z)\coloneqq\int_0^{\infty} \e^{-t}\, t^{z-1}\,\d t\ .
\end{align}
It has the following recursion relation \cite[(5.5.1)]{NIST23}:
\begin{align}\label{eq:gamma func rec}
    \Gamma(z+1)=z\,\Gamma(z)\ .
\end{align}
As special cases, we have:
\begin{align}\label{eq:gamma func 1/2}
    \Gamma(\tfrac{1}{2})=\sqrt{\pi}\ ,\qquad \Gamma(1)=1\ .
\end{align}
When $\Re\, z\leq 0$, $\Gamma(z)$ is defined by analytic continuation through the reflection formula \cite[(5.5.3)]{NIST23}:
\begin{align}\label{eq:gamma func reflection formula}
    \Gamma(z)\,\Gamma(1-z)=\frac{\pi}{\sin(\pi z)}\ .
\end{align}
As clear from \eqref{eq:gamma func reflection formula}, $\Gamma(z)$ has simple poles of residue $(-1)^n/n!$ at $z=-n$ with $n=0,1,2,\cdots$. We further record two relations for the gamma functions relevant in this paper. One is the Euler's beta integral \cite[(5.12.1)]{NIST23}:
\begin{align}\label{eq:Euler's beta integral} 
    \int_0^1\, t^{x-1}(1-t)^{y-1}\,\d t=\frac{\Gamma(x)\,\Gamma(y)}{\Gamma(x+y)}\ ,\qquad \Re\, x>0 \ ,\Re\, y>0\ .
\end{align}
The other is the duplication formula \cite[(5.5.5)]{NIST23}:
\begin{align}\label{eq:Lenegdre duplication formula}
        \Gamma(2z)=\frac{2^{2z-1}}{\sqrt{\pi}}\,\Gamma(z)\,\Gamma(z+\tfrac{1}{2})\ ,\qquad 2z\not\in\{ 0,-1,-2,\cdots\} \ .
\end{align}

We will use standard notations for generalized hypergeometric functions \cite[(16.2.1)]{NIST23}:
\begin{align}\label{eq:generalized hypergeometric}
{}_rF_s\left({a_1,\cdots, a_r\atop b_1,\cdots, b_s};z\right)\coloneqq\sum_{n=0}^{\infty}\,\frac{(a_1)_n\cdots (a_r)_n}{(b_1)_n\cdots (b_s)_n}\,\frac{z^n}{n!}\ ,
\end{align}
for $b_1,\cdots,b_s\not \in \{0,-1,-2,\cdots\}$ with $(a)_n$ being Pochhammer’s symbol defined by \cite[(5.2.4) and (5.2.5)]{NIST23}:
\begin{align}
    (a)_n\coloneqq 
    \begin{dcases}
        a\,(a+1)\cdots (a+n-1)&\\
     \frac{\Gamma(a+n)}{\Gamma(a)}   & a\neq0,-1,-2,\cdots
    \end{dcases}\ .
\end{align}
It is clear from the definition \eqref{eq:generalized hypergeometric} that we have the following confluence property:
\begin{align}\label{eq:confluence property}
    {}_{r+1}F_{s+1}\left({a_1,\cdots, a_r,c\atop b_1,\cdots, b_s,c};z\right)={}_rF_s\left({a_1,\cdots, a_r\atop b_1,\cdots, b_s};z\right)\ .
\end{align}
It is often the case that generalized hypergeometric functions \eqref{eq:generalized hypergeometric} reduce to products of gamma functions when $r=s+1$ and $z=1$. One of the prominent examples is the Pfaff–Saalschütz balanced sum \cite[(16.4.3)]{NIST23}:
\begin{align}\label{eq:balanced sum}
    {}_3F_2\left({-n,a,b\atop c,d};1\right)=\frac{(c-a)_n(c-b)_n}{(c)_n(c-a-b)_n}\ ,
\end{align}
where $c+d=a+b+1-n$ and $n=0,1,2,\cdots$.
We particularly use the special notation for the generalized hypergeometric functions \eqref{eq:generalized hypergeometric} when $r=s=1$ and call it Kummer's confluent hypergeometric function of the first kind \cite[(13.2.2)]{NIST23}:
\begin{align}\label{eq:kummer confluent M def}
    M(a,b,z)\coloneqq {}_1F_1\left({a\atop b};z\right)=\sum_{n=0}^{\infty}\,\frac{(a)_n}{(b)_n}\,\frac{z^n}{n!} \ .
\end{align}
It has the following integral representation \cite[(13.4.1)]{NIST23}:
\begin{align}\label{eq:confluent M integral repr}
    M(a,b,z)=\frac{\Gamma(b)}{\Gamma(a)\,\Gamma(b-a)}\,\int_0^1 \e^{z t}\, t^{a-1}(1-t)^{b-a-1}\,\d t \ ,
\end{align}
for $0<\Re\,a<\Re\,b$.

The Whittaker-$M$ confluent hypergeometric function has the expression in terms of Kummer's confluent hypergeometric function of the first kind \cite[(13.14.2)]{NIST23}:
\begin{align}\label{eq:Whittaker function definition}
     M_{\kappa,\mu}(z)\coloneqq \e^{-\frac{1}{2}z}\,z^{\mu+\frac{1}{2}}\, M(\tfrac{1}{2}+\mu-\kappa,2\mu+1,z)\ .
\end{align}
It satisfies the Whittaker differential equation \cite[(13.14.1)]{NIST23}:
\begin{align}\label{eq:Whittaker def eq}
    \frac{\d^2 W}{\d z^2}+\left(-\frac{1}{4}+\frac{\kappa}{z}+\frac{\frac{1}{4}-\mu^2}{z^2}\right)\,W=0\  ,
\end{align}
with the boundary condition:
\begin{align}\label{eq:Whittaker M boundary condition}
    M_{\kappa,\mu}(z)\xrightarrow[z\to0]{}   z^{\mu+1/2} \ .
\end{align}
One can express various special functions as special cases of the Whittaker-$M$ functions \cite[section 13.18]{NIST23}, such as the lower incomplete gamma function \cite[(8.2.1)]{NIST23}, the error function \cite[(7.2.1)]{NIST23}, and the modified Bessel function of the first kind \cite[(10.2.2)]{NIST23}, each defined by:
\begin{align}
    \gamma(a,z)&\coloneqq \int_0^z t^{a-1}\,\e^{-t} \,\d t \ ,\qquad \Re\,a>0\ , \label{eq:lower incomplete gamma def}\\
    \mathrm{erf}(z)& \coloneqq \frac{2}{\sqrt{\pi}}\,\int_0^z \e^{-t} \,\d t\ , \label{eq:error function def}
 \end{align}
 and 
 \begin{align}\label{eq:modified bessel first def}
      I_\nu(z)& \coloneqq  \frac{(\tfrac{1}{2}z)^\nu}{\Gamma(\nu+1)}\, {}_0F_1\left({-\atop \nu+1};\frac{1}{4}z^2\right)=\sum_{n=0}^{\infty}\,\frac{(\tfrac{1}{2}z)^{\nu+2n}}{\Gamma(\nu+n+1)\, n!} \ .
\end{align}
According to equations (13.18.4), (13.18.6) and (13.18.8) in \cite{NIST23}, one has:
\begin{align}
M_{\mu-\frac{1}{2},\mu}(z)&=2\mu\,\e^{\frac{1}{2}z}\,z^{\frac{1}{2}-\mu}\,\gamma(2\mu,z)\ , \label{eq:lower incomplete gamma ref}\\
M_{-\frac{1}{4},\frac{1}{4}}(z^2)&=\frac{1}{2}\,\e^{\frac{1}{2}z^2}\,\sqrt{\pi z}\,\mathrm{erf}(z)\ , \label{eq:error function rel}\\
M_{0,\nu}(2z)&=2^{2\nu+\frac{1}{2}}\,\Gamma(1+\nu)\,\sqrt{z}\,I_{\nu}(z) \label{eq:modified bessel first rel}\ .
\end{align}
The Whittaker-$M$ function is also related to the hyperbolic sine function by:
\begin{align}\label{eq:hyperbolic sine rel}
    M_{0,\frac{1}{2}}(2z)&=2\sinh z \ .
\end{align}

We now introduce the Bessel function of the first kind by \cite[(10.16.9) and (10.2.2)]{NIST23}:
\begin{align}\label{eq:Bessel J series}
        J_{\nu}(z)\coloneqq \frac{(\tfrac{1}{2}z)^\nu}{\Gamma(\nu+1)}\, {}_0F_1\left({-\atop \nu+1};-\frac{1}{4}z^2\right)=\sum_{n=0}^{\infty}\,\frac{(-1)^n\,(\tfrac{1}{2}z)^{\nu+2n}}{\Gamma(\nu+n+1)\, n!}\ ,
\end{align}
which is subject to the Bessel differential equation \cite[(10.2.1)]{NIST23}:
\begin{align}\label{eq:Besse dif eq}
   z^2\,\frac{d^2 w}{\d z^2}+z\,\frac{\d w}{\d z}+(z^2-\nu^2)\,w=0\ .
\end{align}
The Bessel function of the first kind has the expression in terms of Kummer's confluent hypergeometric function of the first kind \cite[(10.16.5)]{NIST23}:
\begin{align}\label{eq:Bessel J from Whittaker M}
J_\nu(z)=\frac{(\tfrac{1}{2}z)^\nu\,\e^{\mp \i z}}{\Gamma(\nu+1)}\, M(\nu+\tfrac{1}{2},2\nu+1,\pm2\i z)
\ .
\end{align}
The Mellin-Barnes type representation of the Bessel function of the first kind reads \cite[(10.9.22)]{NIST23}:
 \begin{align}\label{eq:Mellin Bessel J}
     J_\nu(x)=\frac{1}{2\pi\i}\,\int_{-\i\,\infty}^{\i\,\infty}\,\frac{\Gamma(-t)\, (\frac{1}{2}x)^{\nu+2t}}{\Gamma(\nu+t+1)}\,\d t\ ,
 \end{align}
for $\Re\,\nu>0,x\in\mathbb{R}_{>0}$.
The Bessel function of the first kind encompasses the trigonometric sine function as a special case \cite[10.16.1]{NIST23}:
\begin{align}\label{eq:Bessel function first sine}
    J_{\frac{1}{2}}(z)=\sqrt{\frac{2}{\pi z}}\,\sin z\ .
\end{align}

This paper proves two identities concerning confluent hypergeometric functions and Bessel functions. They are fundamental and should have been discovered over one hundred years ago by anyone else. However, we could not find either of these identities anywhere in literature, such as \cite{NIST23} and \cite{Watson95}.

\section{Integral representation of the Whittaker-$M$ confluent hypergeometric function}\label{app:Whittaker new integral repr}
This section is devoted to a proof of the following identity:
\begin{theorem}\label{th:Whittaker new integral repr}
    The following integral representation for the Whittaker-$M$ confluent hypergeometric function:
    \begin{align}\label{eq:Whittaker new integral repr}
\begin{aligned}
     M_{\kappa,\mu}(z)&=\frac{\sqrt{\pi }\,\Gamma(2\mu+1)}{2^{\mu}\,\Gamma(\frac{1}{2}(\mu+\kappa+\tfrac{1}{2}))\,\Gamma(\frac{1}{2}(\mu-\kappa+\tfrac{1}{2}))}\\
     &\qquad\times\sqrt{z}\,\int_0^1\xi^{\frac{-\kappa+1/2}{2}-1}(1-\xi)^{\frac{\kappa+1/2}{2}-1}\, \e^{(\xi-1/2) z}\,J_{\mu}(\sqrt{\xi(1-\xi)}z)\,\d\xi\ ,
\end{aligned}
\end{align}
holds true for $\Re\,(\mu\pm\kappa+1/2)>0,z\in\mathbb{C}$.
\end{theorem}
Notice that all arguments of the gamma functions in \eqref{eq:Whittaker new integral repr}, including the one in the numerator $\Gamma(2\mu+1)$, have positive real parts due to the conditions $\Re\,(\mu\pm\kappa+\tfrac{1}{2})>0$, so that one has $\mu\neq\{-\frac{1}{2},-1,-\frac{3}{2},\cdots\}$.
\begin{proof}[Proof of Theorem \ref{th:Whittaker new integral repr}]
We validate the identity by showing that the right-hand side of \eqref{eq:Whittaker new integral repr} satisfies the Whittaker differential equation \eqref{eq:Whittaker def eq} and the boundary condition \eqref{eq:Whittaker M boundary condition}.    
\paragraph{Boundary condition.}
The leading term of the right-hand side of \eqref{eq:Whittaker new integral repr} in taking $z\to0$ limit is:
\begin{align}
& \frac{\sqrt{\pi }\,z^{\mu+\frac{1}{2}}\,\Gamma(2\mu+1)}{2^{2\mu}\,\Gamma(\frac{1}{2}(\mu+\kappa+\tfrac{1}{2}))\,\Gamma(\frac{1}{2}(\mu-\kappa+\tfrac{1}{2}))\,\Gamma(\mu+1)}\,\int_0^1\xi^{\frac{-\kappa+\mu+1/2}{2}-1}(1-\xi)^{\frac{\kappa+\mu+1/2}{2}-1}\,\d\xi\ ,
\end{align}
where we used \eqref{eq:Bessel J series}. Given the condition $\Re\,(\mu\pm\kappa+\tfrac{1}{2})>0$, one can perform the integral concerning $\xi$ by use of \eqref{eq:Euler's beta integral} to have:
\begin{align}
\frac{\sqrt{\pi }\,z^{\mu+\frac{1}{2}}\,\Gamma(2\mu+1)}{2^{2\mu}\,\Gamma(\nu+1)\,\Gamma(\mu+\tfrac{1}{2})} \ .
\end{align}
We can simplify this expression with the aid of \eqref{eq:Lenegdre duplication formula} to arrive at:
\begin{align}
z^{\mu+\frac{1}{2}}\ .
\end{align}
Hence, both sides of \eqref{eq:Whittaker new integral repr} have the same behavior in the limit $z\to0$.

\paragraph{Differential equation.}
Let us check that the right-hand side of \eqref{eq:Whittaker new integral repr} satisfies the Whittaker differential equation \eqref{eq:Whittaker def eq}. By acting the Whittaker differential operator:
\begin{align}
    \frac{\d^2 }{\d z^2}-\frac{1}{4}+\frac{\kappa}{z}+\frac{\frac{1}{4}-\mu^2}{z^2}\  ,
\end{align}
on the right-hand side of \eqref{eq:Whittaker new integral repr} and making some manipulations, one finds that the result is proportional to the following factor:
\begin{align}\label{eq:diff eq check 1}
    \begin{aligned}
&\int_0^1\xi^{\frac{-\kappa+1/2}{2}-1}(1-\xi)^{\frac{\kappa+1/2}{2}-1}\, \e^{(\xi-1/2) z}\\
        &\quad \times\left(\frac{\d^2 }{\d z^2}+\frac{1}{z}\frac{\d}{\d z}+\xi(1-\xi)-\frac{\mu^2}{z^2}+(2\xi-1)\,\frac{\d}{\d z}-2\xi(1-\xi)+\frac{\kappa+\xi-\tfrac{1}{2}}{z}\right)\\
        &\qquad\qquad\qquad\qquad\qquad\qquad\qquad\qquad\qquad\qquad \times J_{\mu}(\sqrt{\xi(1-\xi)}z)\,\d\xi\ .
    \end{aligned}
\end{align}
Namely, the right-hand side of \eqref{eq:Whittaker new integral repr} fulfills the Whittaker differential equation \eqref{eq:Whittaker def eq} if and only if the above expression \eqref{eq:diff eq check 1} vanishes.
Be aware that the first four terms of \eqref{eq:diff eq check 1} add up to zero thanks to the Bessel differential equation \eqref{eq:Besse dif eq}. We now focus on the fifth term. After making use of the identity:
\begin{align}
        (2\xi-1)\,\frac{\d}{\d z}\,J_{\mu}(\sqrt{\xi(1-\xi)}z)&=-\frac{2\xi(1-\xi)}{z}\,\frac{\d}{\d \xi}\,J_{\mu}(\sqrt{\xi(1-\xi)}z)\ ,
\end{align}
that follows immediately from \eqref{eq:Bessel J series}, one finds that the expression \eqref{eq:diff eq check 1} boils down to:
\begin{align}\label{eq:diff eq check 2}
    \begin{aligned}
&\int_0^1\xi^{\frac{-\kappa+1/2}{2}-1}(1-\xi)^{\frac{\kappa+1/2}{2}-1}\, \e^{(\xi-1/2) z}\\
        &\qquad\times\left(-\frac{2\xi(1-\xi)}{z}\,\frac{\d}{\d \xi}\,-2\xi(1-\xi)+\frac{\kappa+\xi-\tfrac{1}{2}}{z}\right)\, J_{\mu}(\sqrt{\xi(1-\xi)}z)\,\d\xi\ .
    \end{aligned}
\end{align}
Let us perform integration by parts for $\xi$ provided the condition $\Re\,(\mu\pm\kappa+\tfrac{1}{2})>0$ necessary to drop off the surface terms. Then, we see that the first term in \eqref{eq:diff eq check 2} turns out to take the form:
\begin{align}
    \begin{aligned}
\int_0^1\xi^{\frac{-\kappa+1/2}{2}-1}&(1-\xi)^{\frac{\kappa+1/2}{2}-1}\, \e^{(\xi-1/2) z}\,\left(2\xi(1-\xi)-\frac{\kappa+\xi-\tfrac{1}{2}}{z}\right)\,J_{\mu}(\sqrt{\xi(1-\xi)}z)\,\d\xi\ .
    \end{aligned}
\end{align}
Substituting this expression for the first term in \eqref{eq:diff eq check 2}, one can readily see that all terms cancel out and conclude that the right-hand side of \eqref{eq:Whittaker new integral repr} satisfies the Whittaker differential equation \eqref{eq:Whittaker def eq}.
This completes the proof of Theorem \ref{th:Whittaker new integral repr}.
\end{proof}

In what follows, we study several special cases of the identity \eqref{eq:Whittaker new integral repr}.
\begin{corollary}\label{th:lower incomplete}
 The lower incomplete gamma function has the following integral representation for $\Re\,\mu>0$:
    \begin{align}\label{eq:lower incomplete integral Bessel}
        \gamma(2\mu,z) =\frac{\Gamma(2\mu+1)\,z^{\mu}\,\e^{-z}}{2^{\mu+1}\,\Gamma(\mu+1)}\,\int_0^1\xi^{\frac{-\mu+1}{2}-1}(1-\xi)^{\frac{1}{2}\mu-1}\, \e^{\xi z}\,J_{\mu}(\sqrt{\xi(1-\xi)}z)\,\d\xi\ .
\end{align}
\end{corollary}

\begin{proof}[Proof of Corollary \ref{th:lower incomplete}]
Owing to Theorem \ref{th:Whittaker new integral repr}, we have for $\Re\,\mu>0$:
\begin{align}
         M_{\mu-\frac{1}{2},\mu}(z)=\frac{\Gamma(2\mu+1)\,\sqrt{z}}{2^{\mu}\,\Gamma(\mu)}\,\int_0^1\xi^{\frac{-\mu+1}{2}-1}(1-\xi)^{\frac{1}{2}\mu-1}\, \e^{(\xi-1/2) z}\,J_{\mu}(\sqrt{\xi(1-\xi)}z)\,\d\xi\ .
\end{align}
We compare this expression with \eqref{eq:lower incomplete gamma ref} using the recursion relation for the gamma function \eqref{eq:gamma func rec} to conclude \eqref{eq:lower incomplete integral Bessel}.
\end{proof}

\begin{corollary}\label{th:error function}
 The error function has the following integral representation:
    \begin{align}\label{eq:error function integral Bessel}
\mathrm{erf}(z) =\frac{ \sqrt{z}\,\e^{-z^2}}{2^{\frac{1}{4}}\,\Gamma(\tfrac{1}{4})} \,\int_0^1\xi^{\frac{3}{8}-1}(1-\xi)^{\frac{1}{8}-1}\, \e^{\xi  z^2}\,J_{\frac{1}{4}}(\sqrt{\xi(1-\xi)}z^2)\,\d\xi    \ .
\end{align}
\end{corollary}

\begin{proof}[Proof of Corollary \ref{th:error function}]
When $\kappa=-\frac{1}{4}$ and $\mu=\tfrac{1}{4}$, the identity \eqref{eq:Whittaker new integral repr} is simplified to give:
\begin{align}
    M_{-\frac{1}{4},\frac{1}{4}}(z^2)=\frac{\sqrt{\pi}\, z\,\e^{-\frac{1}{2}z^2}}{2^{\frac{5}{4}}\,\Gamma(\tfrac{1}{4})} \,\int_0^1\xi^{\frac{3}{8}-1}(1-\xi)^{\frac{1}{8}-1}\, \e^{\xi  z^2}\,J_{\frac{1}{4}}(\sqrt{\xi(1-\xi)}z^2)\,\d\xi\ .
\end{align}
where we used \eqref{eq:gamma func rec} and \eqref{eq:gamma func 1/2} to have $\Gamma(\tfrac{3}{2})=\tfrac{1}{2}\,\Gamma(\tfrac{1}{2})=\frac{1}{2}\,\sqrt{\pi}$. 
By substituting this expression into \eqref{eq:error function rel}, we end up with \eqref{eq:error function integral Bessel}.
\end{proof}

\begin{corollary}\label{th:modified Bessel}
 The modified Bessel function of the first kind has the following integral representation for $\Re\,\nu>-\tfrac{1}{2}$:
    \begin{align}\label{eq:modified Bessel integral Bessel}
   I_{\nu}(z)=  \frac{\Gamma(\tfrac{1}{2}(\nu+\tfrac{3}{2}))}{\sqrt{2\pi}\,\Gamma(\tfrac{1}{2}(\nu+\tfrac{1}{2}))}\,\int_0^1\xi^{\frac{1}{4}-1}(1-\xi)^{\frac{1}{4}-1}\, \e^{(2\xi-1) z}\,J_{\nu}(2z\sqrt{\xi(1-\xi)})\,\d\xi\ .
\end{align}
\end{corollary}

\begin{proof}[Proof of Corollary \ref{th:modified Bessel}]
It follows from Theorem \ref{th:Whittaker new integral repr} and \eqref{eq:modified bessel first def} that:
\begin{align}
\begin{aligned}
       I_{\nu}(z)&=
   \frac{\sqrt{\pi }\,\Gamma(2\nu+1)}{2^{3\nu}\,\Gamma(\nu+1)\,\Gamma(\frac{1}{2}(\nu+\tfrac{1}{2}))\,\Gamma(\frac{1}{2}(\nu+\tfrac{1}{2}))}\\
     &\qquad\times\int_0^1\xi^{\frac{1}{4}-1}(1-\xi)^{\frac{1}{4}-1}\, \e^{(2\xi-1) z}\,J_{\nu}(2z\sqrt{\xi(1-\xi)})\,\d\xi\ ,
\end{aligned}
\end{align}
for $\Re\,\nu>-\tfrac{1}{2}$.
The duplication formula \eqref{eq:Lenegdre duplication formula} allows us to simplify the gamma factors in the above expression, resulting in \eqref{eq:modified Bessel integral Bessel}.
\end{proof}

\begin{corollary}\label{th:hyperbolic sine}
 The hyperbolic sine function has the following integral representation:
    \begin{align}\label{eq:hyperbolic sine integral Bessel}
\sinh z=\frac{1}{2\pi}\int_0^1 \frac{\e^{(2\xi-1)z}\,\sin(2\sqrt{\xi(1-\xi)}z)}{\xi(1-\xi)}\,\d\xi\ .
\end{align}
\end{corollary}

\begin{proof}[Proof of Corollary \ref{th:hyperbolic sine}]
This identity follows immediately from \eqref{eq:Whittaker new integral repr} by setting $\kappa=0,\mu=\tfrac{1}{2}$ and using \eqref{eq:Bessel function first sine} and \eqref{eq:hyperbolic sine rel}.
\end{proof}

\section{Summation formula for Kummer's confluent hypergeometric functions of the first kind}\label{app:Whittaker new sum repr}
In this section, we first roughly manipulate the identity \eqref{eq:Whittaker new integral repr} to find that it implies a duplication-type summation formula for Kummer's confluent hypergeometric function of the first kind for a limited parameter region. We then move on to a rigorous proof of the identity for more generic arguments. Along the way, we substantiate a finite summation formula involving Pochhammer symbols. 

Assuming that $\Re\,(\mu\pm\kappa+\tfrac{1}{2}),\Re\,\mu>0$ and $z\in\mathbb{R}_{>0}$, we apply the Mellin-Barnes integral representation of the Bessel function of the first kind \eqref{eq:Mellin Bessel J} to \eqref{eq:Whittaker new integral repr}.
After the change of the order of integration, we find that the second line of \eqref{eq:Whittaker new integral repr} takes the form:
\begin{align}\label{eq:rough mani after mellin 1}
\frac{1}{2\pi\i}\,\int_{-\i\,\infty}^{\i\,\infty}\,\frac{\Gamma(-t)\,\sqrt{z}\, \e^{-\frac{1}{2}z}\, (\tfrac{1}{2}z)^{\mu+2 t}}{\Gamma(\mu+t+1)} \,\int_0^1\xi^{\frac{\mu-\kappa+1/2}{2}+t-1}(1-\xi)^{\frac{\mu+\kappa+1/2}{2}+t-1}\, e^{z \xi}\,\d\xi\,\d t\ .
 \end{align}
Notice that one can identify the $\xi$-integral in the above expression with the integral representation of Kummer's confluent hypergeometric functions of the first kind \eqref{eq:confluent M integral repr}. Hence, after some calculations using the duplication formula \eqref{eq:Lenegdre duplication formula}, we see that the expression \eqref{eq:rough mani after mellin 1} becomes:
\begin{align}
 \begin{aligned}
\frac{2^\mu\,z^{\mu+\frac{1}{2}}\,\e^{-\frac{1}{2}z}}{\sqrt{\pi}}\,\frac{1}{2\pi\i}\, \int_{-\i\,\infty}^{\i\,\infty}&\frac{\Gamma(\tfrac{\mu-\kappa+1/2}{2}+t)\,\Gamma(\tfrac{\mu+\kappa+1/2}{2}+t)\,\Gamma(\mu+\tfrac{1}{2}+t)\,\Gamma(-t)}{\Gamma(2\mu+1+2t)\,\Gamma(\mu+\frac{1}{2}+2t)}\\
      &\qquad\qquad\times  z^{2 t}\,M(\tfrac{\mu-\kappa+1/2}{2}+t,\mu+\tfrac{1}{2}+2t,z)\,\d t\ .
 \end{aligned}
 \end{align}
We now deform the integration contour to the right and pick up the residues from $\Gamma(-t)$. Then, the second line of \eqref{eq:Whittaker new integral repr} turns out to be:
\begin{align}\label{eq:Whittaker new sum repr 0}
 \begin{aligned}
    \frac{2^\mu\, z^{\mu+\frac{1}{2}}\,e^{-\frac{1}{2}z}}{\sqrt{\pi}}\,&\sum_{n=0}^{\infty}\,\frac{\Gamma(\tfrac{\mu-\kappa+1/2}{2}+n)\,\Gamma(\tfrac{\mu+\kappa+1/2}{2}+n)\,(\mu+\tfrac{1}{2})_{n}}{(\mu+\tfrac{1}{2})_{2n}\,\Gamma(2\mu+1+2n)\, n!}\\
    &\qquad\qquad \times (-z^2)^{n}\,M(\tfrac{\mu-\kappa+1/2}{2}+n,\mu+\tfrac{1}{2}+2n,z)\ .
 \end{aligned}
 \end{align}
 We substitute this expression for the second line of \eqref{eq:Whittaker new integral repr} and compare both sides using \eqref{eq:Whittaker function definition} to conclude that:
\begin{align}\label{eq:Whittaker new sum repr 1}
     M(2a,2b,z)=\sum_{n=0}^{\infty}\,\frac{(a)_n(b)_n(b-a)_n}{(b)_{2n}(2b)_{2n}\, n!} \, (-z^2)^{n}\,M(a+n,b+2n,z) \ ,
\end{align}  
  for $z\in\mathbb{R}_{>0}$ and $\Re\,(a),\Re\,(b-a),\Re\,(b+\tfrac{1}{2})>0$. 

 We can verify that \eqref{eq:Whittaker new sum repr 1} holds true for any $a,b,z\in\mathbb{C}$. To this end, we start by showing the following proposition:
 \begin{proposition}\label{pro:identity}
 The following finite summation formula involving Pochhammer symbols holds true for $k\in\{0,1,2,\cdots\}$:
  \begin{align}\label{eq:identity}
 \frac{(2a)_k}{(2b)_k}=\frac{k!}{(b)_k}\,\sum_{n=0}^{\lfloor \frac{1}{2}k\rfloor}\,\frac{(a)_{k-n}\, (b)_n\,(b-a)_n}{(2b)_{2n}\, n!\, (k-2n)!}\ ,
 \end{align}     
where $\lfloor x\rfloor$ is the floor function that returns the largest integer less than or equal to $x$.
 \end{proposition}

\begin{proof}[Proof of Proposition \ref{pro:identity}]
 First of all, implementing the duplication formula \eqref{eq:Lenegdre duplication formula} and reflection formula \eqref{eq:gamma func reflection formula}, we find that the right-hand side of \eqref{eq:identity} boils down to a generalized hypergeometric function with unit argument:
 \begin{align}\label{eq:identity 1}
\frac{(a)_k}{(b)_k}\,{}_3F_2\left({b-a,-\tfrac{1}{2}k+\tfrac{1}{2},-\tfrac{1}{2}k\atop b+\frac{1}{2},1-a-k};1\right)\ .
\end{align}
Take care that this expression is a finite sum, as either of the top parameters $-\tfrac{1}{2}k+\tfrac{1}{2}$ and $-\tfrac{1}{2}k$ is always a non-positive integer and the series expansion terminates. One can further reduce the expression \eqref{eq:identity 1} utilizing Pfaff–Saalschütz balanced sum \eqref{eq:balanced sum} as shown below.
When $\tfrac{1}{2}k\in\{0,1,2,\cdots\}$, owing to \eqref{eq:balanced sum}, the expression \eqref{eq:identity 1} reduces to:
\begin{align}
    \frac{(a)_k\,(a+\tfrac{1}{2})_{\frac{1}{2}k}\,(b+\tfrac{1}{2}k)_{\frac{1}{2}k}}{(b)_k\,(b+\tfrac{1}{2})_{\frac{1}{2}k}\,(a+\tfrac{1}{2}k)_{\frac{1}{2}k}}\ .
\end{align}
Through the duplication formula \eqref{eq:Lenegdre duplication formula}, we can simplify this expression to give the same one as the left-hand side of \eqref{eq:identity}.
One can perform a similar analysis for the case $\tfrac{1}{2}k-\tfrac{1}{2}\in\{0,1,2,\cdots\}$ and verify that \eqref{eq:identity} is an identity for $k\in\{0,1,2,\cdots\}$.
\end{proof}

We are now in a position to prove the following:
\begin{theorem}\label{th:Whittaker new sum repr}
    The following duplication-type summation formula for Kummer's confluent hypergeometric function of the first kind holds true for $a,b,z\in\mathbb{C}$:
\begin{align}\label{eq:Whittaker new sum repr}
     M(2a,2b,z)=\sum_{n=0}^{\infty}\,\frac{(a)_n(b)_n(b-a)_n}{(b)_{2n}(2b)_{2n}\, n!} \, (-z^2)^{n}\,M(a+n,b+2n,z) \ .
\end{align}
\end{theorem}
This summation formula follows immediately from Proposition \ref{pro:identity} as below.
\begin{proof}[Proof of Theorem \ref{th:Whittaker new sum repr}]
    By expanding Kummer's confluent hypergeometric functions of the first kind \eqref{eq:kummer confluent M def} in power series and making a change of summation variables, one finds that the right-hand side of \eqref{eq:Whittaker new sum repr} becomes:
    \begin{align}
        \sum_{k=0}^{\infty}\,\frac{z^k}{(b)_k}\,\sum_{n=0}^{\lfloor \frac{1}{2}k\rfloor}\,\frac{(a)_{k-n}\, (b)_n\,(b-a)_n}{(2b)_{2n}\, n!\, (k-2n)!}\ .
    \end{align}
    Using the result of Proposition \ref{pro:identity}, this expression reduces to:
     \begin{align}
        \sum_{k=0}^{\infty}\,\frac{(2a)_k}{(2b)_k}\,\frac{z^k}{k!}\ .
    \end{align}
Bearing in mind the definition of Kummer's confluent hypergeometric function of the first kind \eqref{eq:kummer confluent M def}, one can equate this to the left-hand side of \eqref{eq:Whittaker new sum repr}.
    This completes the proof of Theorem \ref{th:Whittaker new sum repr}.
\end{proof}

The following summation formula for the Bessel function of the first kind follows immediately from Theorem \ref{th:Whittaker new sum repr}:
\begin{corollary}\label{th:Bessel J new sum rule}
 The following summation formula for the Bessel function of the first kind holds true:
    \begin{align}\label{eq:Bessel J new sum rule}
    J_{2\nu+\frac{1}{2}}(z)=\frac{ \Gamma(\nu+1)}{\Gamma(2\nu+\tfrac{3}{2})}\,\sum_{n=0}^{\infty}\,\frac{(\nu+\tfrac{1}{2})_n}{(2\nu+\tfrac{3}{2})_n\, n!} \, \left(\frac{z}{2}\right)^{\nu+\frac{1}{2}+n}\,J_{\nu+n}(z)\ ,
\end{align}
for $\nu\not\in\{-1,-2,\cdots\}$.
\end{corollary}

We notify that $\nu$ must not be negative integers where the numerator gamma function $\Gamma(\nu+1)$ on the right-hand side of \eqref{eq:Bessel J new sum rule} would not be defined.
We further remark that owing to \eqref{eq:Bessel function first sine}, when $\nu=0$ the summation formula \eqref{eq:Bessel J new sum rule} reduces to the known expansion of the trigonometric sine function in terms of the Bessel function of the first kind \cite[(9.4.2.19)]{Luke69}:
\begin{align}\label{eq:sin expansion Bessel J}
\frac{\sin z}{z} = \sum_{n=0}^{\infty}\,\frac{(\tfrac{1}{2}z)^{n}}{n!\,(2n+1)}\,J_{n}(z)\ .
\end{align}

\begin{proof}[Proof of Corollary \ref{th:Bessel J new sum rule}]
    It follows from the relation \eqref{eq:Bessel J from Whittaker M} that the left-hand side of \eqref{eq:Bessel J new sum rule} becomes:
    \begin{align}
 \frac{\e^{\mp \i z}}{\Gamma(2\nu+\tfrac{3}{2})} \,\left(\frac{z}{2}\right)^{2\nu+\frac{1}{2}}\, M(2\nu+1,4\nu+2,\pm 2\i z)\ .
\end{align}
Employing \eqref{eq:Whittaker new sum repr} and \eqref{eq:Bessel J from Whittaker M}, one sees that this expression becomes:
\begin{align}
    \begin{aligned}
     &\frac{ \Gamma(\nu+1)}{\Gamma(2\nu+\tfrac{3}{2})}\,\sum_{n=0}^{\infty}\,\frac{2^{4n}\,(\nu+\tfrac{1}{2})_n^2(\nu+1)_n(2\nu+1)_n}{(2\nu+1)_{2n}(4\nu+2)_{2n}\, n!} \,\left(\frac{z}{2}\right)^{\nu+\frac{1}{2}+n}\,J_{\nu+n}(z)\ .
    \end{aligned}
\end{align}
We can simplify this expression by applying the duplication formula \eqref{eq:Lenegdre duplication formula} several times. We then finally arrive at \eqref{eq:Bessel J new sum rule}.
\end{proof}

We here make clear that Corollary \ref{th:Bessel J new sum rule} also follows from the identity found in \cite[(9.4.7.2)]{Luke69}:\footnote{We are grateful to Roy Hughes for pointing out that our result \eqref{eq:sin expansion Bessel J} is included in this known identity.}
\begin{align}\label{eq:2F1 Bessel J}
    {}_1F_2\left({a\atop b,c};-\frac{1}{4}z^2\right)=\Gamma(c)\,\left(\frac{2}{z}\right)^{c-1}\,\sum_{n=0}^{\infty}\,\frac{(b-a)_n}{n!\,(b)_n}\,\left(\frac{z}{2}\right)^n\,J_{n+c-1}(z)\ ,
\end{align}
by setting $a=c=\nu+1,b=2\nu+\tfrac{3}{2}$ with the definition of the Bessel function \eqref{eq:Bessel J series} and the confluence property of generalized hypergeometric functions \eqref{eq:confluence property} in mind. Alternatively but more directly, one can also prove Corollary \ref{th:Bessel J new sum rule} from the identity \cite[(9.4.4.3)]{Luke69}, that is essentially subsumed by \eqref{eq:2F1 Bessel J}:
\begin{align}
    J_{\nu}(z)=\frac{\Gamma(a)}{\Gamma(\nu+1)}\,\left(\frac{z}{2}\right)^{\nu+1-a}\,\sum_{n=0}^{\infty}\,\frac{(\nu+1-a)_n}{n!\,(\nu+1)_n}\,\left(\frac{z}{2}\right)^{n}\, J_{n+a-1}(z)\ ,
\end{align}
with $a\not\in\{0,-1,-2,\cdots\}$ by making the following replacements:
\begin{align}
    \begin{dcases}
        \nu\to 2\nu+\tfrac{1}{2}\\
        a\to \nu+1
    \end{dcases}\ .
\end{align}
All these confirm the overall validity of our findings.

\section{Conclusion}
In this paper, we proved two novel relations concerning confluent hypergeometric functions and Bessel functions (Theorem \ref{th:Whittaker new integral repr} and \ref{th:Whittaker new sum repr}).
We specified their parameters and obtained various identities. To the best of our knowledge, four integration formulas (Corollary \ref{th:lower incomplete}, \ref{th:error function}, \ref{th:modified Bessel} and \ref{th:hyperbolic sine}) that follow directly from Theorem \ref{th:Whittaker new integral repr} are new. On the other hand, Corollary \ref{th:Bessel J new sum rule} is not in literature but is implicit within known formulas. In deriving Theorem \ref{th:Whittaker new sum repr}, we have also shown a finite summation formula for products of Pochhammer symbols (Proposition \ref{pro:identity}). It will be interesting to consider applications or extensions of our findings.

\section{Acknowledgments}
We are grateful to Tatsuma Nishioka for the valuable discussions. We would like to thank Roy Hughes and Kohei Fukai for correspondence and comments on the manuscript. The author is supported by JSPS KAKENHI Grant-in-Aid for JSPS fellows Grant No.\,21J20750, JSR Fellowship, and FoPM (WINGS Program), the University of Tokyo.

\end{document}